\documentclass[12pt]{article}

\usepackage{latexsym,amsfonts,amssymb,theorem,dina4}
\newcommand{\ebox}{\hspace*{\fill{$\Box$}}}
\newcommand{\lr}{\longrightarrow}

\title{A note on the Thom isomorphism in geometric (co)homology}

\newtheorem{thm}{Theorem}[section]

\author{Martin Jakob}

\begin{document}

\maketitle

\begin{abstract}
Using geometric homology and cohomology 
we give a simple conceptual proof of 
the Thom isomorphism theorem. 
\end{abstract}

\section{Introduction}

In the articles \cite{althomol} and 
\cite{geomhomol} we constructed homology functors in a geometric way. 
We showed that these functors which are defined on pairs of spaces
satisfy the Eilenberg--Steenrod axioms for a generalized homology theory. 

Let $h^*$ be a multiplicative cohomology theory (for example 
singular cohomology with coefficients in a ring, or $K$--theory). 
In this setting, the associated homology groups 
$h_q (X)$ of a space $X$ are
defined by means of 
 of triples 
$(M, x, f)$ (so--called {\em cycles}), where $M$ is a manifold, 
$x \in h^* (M)$ is a cohomology class
of $M$ and $f: M \to X$ is a continuous map. 

An equivalence relation on the set 
of cycles must be imposed. It  includes bordism and 
a procedure called ``vector bundle modification''
to shift the degree of $x$ and the dimension of $M$
without leaving the class of $(M, x, f)$. 

In the case of $K$--theory this construction goes back to some work
of Paul Baum, cf.~for example \cite{bd}. 

\bigskip

On the category of differentiable manifolds and smooth maps
 there is also 
a bivariant version of this construction
(cf.~\cite{biv}). 
For an introductory text to bivariant theories see 
\cite{fm}.

\bigskip

These geometric approaches to homology as opposed
to the usual spectral methods of stable homotopy allow
surprisingly simple proofs of 
properties requiring orientations (say of manifolds or maps). 
In the usual stable homotopy setting often extra considerations
are needed. An illustration of this point is the proof
of the Poincar\'e duality theorem.

\bigskip

Recently, the question has been raised whether in the geometric
(co)homology setting  there is a simple proof
of the Thom isomorphism theorem. That is in fact the case and 
in this little note we shall prove: 
\begin{thm}
Let $h^*$ be a multiplicative cohomology theory and 
let $\pi: E \to X$ be a smooth $h^*$--oriented vector bundle of 
rank $n$. Then the geometric Thom class 
$$
t_E = [X, 1_X, \sigma_0] \in h^n (DE, SE)
$$
induces isomorphisms
$$
h^q (X) \to h^{q+n} (DE, SE), \ 
x \mapsto t_E \cdot \pi^* (x) 
$$
 and 
$$
h^{q+n} (DE, SE) \to h^q(X), \ y \mapsto \pi_! (y)
$$
which are inverse to each other. 
\end{thm}
Here $DE$, resp.~$SE$ is the total space of 
the unit ball, resp.~unit sphere 
bundle associated to some metric on $E$, and $\sigma_0$ denotes
the zero section
 of the 
ball bundle. 
Further, the map $\pi_!$ sends the 
geometric class $[P, x, g] \in h^{q+n}$ onto the class 
$[P, x, \pi g] \in h^q(X)$. 

To be precise, in geometric {\em co}homology 
$DE$ should be a manifold without boundary. 
So one has to put $DE =$ open ball bundle in $E$
 of a fixed radius $> 1$. 
 
\section{The proof of the Thom isomorphism theorem}

Let us observe first that $\sigma_0$ in the cycle 
$(X, 1_x, \sigma_0)$ avoids the sphere bundle $SE$ and thus
describes in fact a geometric cohomology class of the pair $(DE, SE)$. 

\bigskip

Let $[M, x, f] \in h^q (X)$. Now
$t_E \cdot \pi^* ([M, x, f]) \in h^{q + n} (DE, SE)$ 
is represented by the composition of pull backs
$$
\begin{array}{ccc}
(M \times_X DE) \times _{DE} X & \buildrel f'' \over \lr & X \\
\downarrow && \downarrow \\
M \times_X DE & \buildrel f' \over \lr & DE \\
\downarrow && \downarrow \\
M & \buildrel f \over \lr & X.
\end{array}
$$

Then 
$$
((M \times_X DE) \times _{DE} X, x \times 1 \times 1, f'')
$$
represents $\pi_! (t_E \cdot \pi^* ([M, x, f])$. The  
equivalence to $(M, x, f)$ is given by the isomorphism
$$
(M \times_X DE) \times _{DE} X \lr M, \ 
((m, v), x) \mapsto m.
$$

Therefore we have shown
$$
\pi_! (t_E \cdot \pi^*( \ldots)) = id_{h^*(X)}.
$$
Let $(P, x, g)$ be a cycle representing a cohomology class $y$ of 
$h^*(DE, SE)$. 
Then $t_E \cdot \pi^*(\pi_!(y))$ is represented by the pull back diagram
$$
\begin{array}{ccc}
DE \times_X P & \lr &  P \\
\downarrow && \downarrow \\
DE & \lr & X \\
\end{array}
$$

Multiplication with the Thom class is given by the pull back 
diagram
$$
\begin{array}{ccc}
X \times_{DE} (DE \times_X P) & \lr & DE \times_X P \\
\downarrow && \downarrow \\
X & \buildrel \sigma_0 \over \lr & DE \\
\end{array}
$$
Now observe that 
$$
X \times _{DE} (DE \times _X P) = 
\{(x, v, p) \in X \times DE \times P; \sigma_0 (x) = v, \pi g (p) = x \}
$$
is diffeomorphic to $P$ by sending $(x, v, p)$ onto $p$. 

\bigskip

To finish the proof, one needs to show
that the map from the last diagram 
$$
X \times _{DE} (DE \times _X P) \lr DE, \ 
(x, (v, p)) \mapsto v
$$
and the map
$$
g: P \lr DE
$$
are cobordant. This is done by a homotopy joining $g (p)$ and
its projection onto the zero section. \ebox

\section{Two final remarks}

1. The proof of the homological Thom isomorphism 
follows a similar pattern. 
Let $[P, x, g] \in h_q (DE, SE)$. Consider the pull back diagram 
$$
\begin{array}{ccc}
P \times_{DE} X & \buildrel g' \over \lr & X \\
\downarrow && \downarrow_{\sigma_0} \\
P &  \buildrel g \over \lr & DE. \\
\end{array}
$$
This gives the cycle $(P \times_{DE} X, x, g')$ representing a class
in $h_{q - n} (X)$. 

\bigskip

On the other hand, for a cycle $(M, x, f)$ representing a 
class in $h_q (X)$ we get a class in $h_{q + n} (DE, SE)$ 
via the pull back 
$$
\begin{array}{ccc}
f^* DE & \lr & DE \\
\downarrow && \downarrow \\
M & \lr & X \\
\end{array}
$$
These two constructions are inverse to each other. 

\bigskip

Our geometric approach to the homological  Thom isomorphism is valid
in a more general setting: Let $E \to X$ be a vector bundle
in the paracompact category. For a class 
$[M, x, f]$ of $h_* (X)$ we can perform the latter pull back, when 
we replace $f: M \to X$ by the composition
$$
M \lr X \lr BO_n, 
$$
where the last arrow is the classifying map of the vector bundle. 

Taking a smooth  universal bundle 
$D EO_n \to BO_n$ we can make the composition
$M \to BO_n$ smooth. Moreover, the pull back is a smooth, finite dimensional
manifold with boundary. In favorable cases it maps to $(DE, SE)$.

\bigskip

2. Clearly, the Thom class $t_E$ can also be viewed as a bivariant class
$$
[X, 1_X, \sigma_0] \in h^n (E \buildrel \pi \over \lr X).
$$
It is not hard to work out the appropriate version of the 
Thom isomorphism theorem using the intersection product of 
geometric bivariant theories.

\bigskip

M. Jakob

Mathematical Institute

Univ.~of Neuchatel

CH--2000 Neuchatel, Switzerland

{\tt martin.jakob$\symbol{64}$unine.ch}


\begin{thebibliography}{1}
\bibitem{bd} Baum, P. and Douglas, R.: {\em $K$ Homology and
Index Theory}, Proceedings of Symposia in Pure Mathematics,
{\bf 38}, Part 1, 117 -- 173 (1982)
\bibitem{fm} Fulton, W. and MacPherson, R.:
{\em Categorical framework for the study of singular spaces}, 
Memoirs of the AMS, {\bf 243}, 1981 
\bibitem{biv} Jakob, M.: {\em Bivariant theories for smooth manifolds}, 
 Appl.~Categ.~Struct.~{\bf 10}, No.3, 279-290 (2002)
\bibitem{althomol} Jakob, M.:
{\em An Alternative Approach to Homology},
Contemp.~Math.~{\bf 265}, 87 -- 97 (2000)
\bibitem{geomhomol} Jakob, M.:
{\em A bordism--type construction of homology},
Manuscr.~Math.~{\bf 96}, 67 -- 80 (1998)
\end{thebibliography}
\end{document}